\documentstyle[12pt]{article}

\font\tenmsb=msbm10 scaled \magstep1
\font\sevenmsb=msbm7 scaled \magstep1
\font\fivemsb=msbm5 scaled \magstep1

\newfam\msbfam
\textfont\msbfam=\tenmsb
\scriptfont\msbfam=\sevenmsb
\scriptscriptfont\msbfam=\fivemsb
\def\Bbb#1{{\fam\msbfam\relax#1}}
\begin{document}

\def\R{{\Bbb R}}
\def\C{{\Bbb C}}
\def\N{{\Bbb N}}
\def\Z{{\Bbb Z}}
\def\Q{{\Bbb Q}}
\def\S{{\Bbb S}^{p-1}}
\def\be{\begin{equation}}
\def\en{\end{equation}}
\def\Subset{\subset\subset}
\def\e{\varepsilon}
\def\d{\delta}
\def\j1{j=1,\dots,}
\def\l{\lambda}
\def\grad{{\rm grad}}
\def\Arg{{\rm Arg}}
\def\Subset{\subset\subset}
\def\sp{{\rm sp}}
\def\mes{{\rm mes}}
\def\apf{f\in AP(T_\Omega)}
\def\OA{\Omega\setminus A_f}
\def\DA{D\setminus A_f}
\def\yD{y\in\overline{D'}}
\def\oQ{\overline{Q}}
\def\L{\Lambda}
\def\tL{\tilde\Lambda}
\def\tm{\tilde m}
\def\D{\Delta}
\def\TD{T_{D'\rho},\ \rho=1,2}
\def\<{\langle}
\def\>{\rangle}
\title{
    \hfill\raisebox{0.0ex}[1.ex][3ex]{}
       {\rm \small  MSC 42A75 (32A60, 32A18)}
      \newline
\bf Holomorphic almost periodic functions in tube domains
and their amoebas
\footnote {This research was
supported by INTAS-99-00089 project}
}

\author{Favorov S.Ju.}

\date{}

\maketitle

\begin{abstract}

We extend the notion of amoeba to holomorphic almost periodic
functions in tube domains. In this setting, the order of
a function in a connected component of the complement to
 its amoeba is just the mean motion of this function.
 We also find  a correlation
between the orders in  different components.\par
\end{abstract}

\medskip


\bigskip


The notion of amoebas for algebraic varieties was introduced
 as auxiliary tools for studying their topological properties.
Then amoebas were studied in various areas of mathematics (algebraic
geometry, topology, combinatorics).
Here we extend this notion to the zero sets of exponential polynomials and,
more generally, to the zero sets of their uniform limits, so-called
holomorphic almost periodic functions. Note that the construction of amoebas
for exponential polynomials is simpler (and more natural) than for usual
 ones. Moreover, it turns out that some notions connected with the amoebas
 coincide in our case with classical notions from the theory of almost
 periodic functions. Hence we believe that our investigation will
be useful both in the theory of holomorphic almost periodic functions and
in the theory of amoebas.

Let $L(w)$ be a Laurent polynomial of $w=(w_1,\dots,w_p)\in\C^p$,
i.e., a finite sum
\be
L(w)=\sum c_mw^m,\quad {\rm where}\ m=(m_1,\dots,m_p)\in\Z^p,
\quad w^m=w_1^{m_1}\dots w_p^{m_p}.
\label{Laur}
\en
The {\it amoeba} $A_L$ of the polynomial $L$ is the image of its zero set
under the map $\alpha: (w_1,\dots,w_p)\mapsto
(\log|w_1|,\dots,\log|w_p|)$ (see \cite{G}). Connected components of
the complement to the amoeba $A_L$ were studied in  \cite{G},
\cite{M}, \cite{P}, \cite{Ru}.
In particular, in the paper \cite{P}
the authors gave an estimate for the number of these components
and introduced the notion of {\it order}
$\nu^k=(\nu_1^k,\dots,\nu_p^k)$ of the  component
$D_k$ by  the formula
\be
\nu_j^k= {1\over 2\pi}\D_{C_j}\Arg L(w),\quad\j1 p,
\label{ord}
\en
where $\D_{C_j}$ means the increment over the circle
$$
C_j=\{w=(w_1^0,\dots,w_{j-1}^0,e^{i\varphi}w_j^0,w_{j+1}^0,\dots,w_p):
\,0\le\varphi\le 2\pi\}
$$
with $w^0=(w_1^0,\dots,w_p^0)\in\alpha^{-1}(D_k)$.
It is easy to see that the numbers $\nu_j^k$ are integers for all
$\j1,p$ and do not depend on $w^0\in\alpha^{-1}(D_k)$.

Observe that after the substitution $w_j=e^{iz_j},\ \j1 p$, the sum
(\ref{Laur}) takes the form
\be
P(z)=\sum c_\l\exp\{i\<z,\l\>\},
\label{sumexp}
\en
with $\l=m\in\Z^p$;  now $D_k$ is a connected component of the set
$$
\{y\in\R^p:\ P(x-iy)\neq 0\quad\forall x\in\R^p\}
$$
and formula (\ref{ord}) takes the form
\be
\nu_j^k={1\over 2\pi}\D_{|x_j|\le\pi}\Arg P(x+iy),
\label{ord1}
\en
where $x=(x_1,\dots,x_p)$ and $-y\in D_k$.

Now remove the condition $\l\in\Z$ and consider  the sum
(\ref{sumexp}) with arbitrary $\l\in\R^p$. More generally,  let
$T_\Omega=\{z=x+iy:\ x\in\R^p,\ y\in \Omega\subset \R^p\}$ be a tube domain
in $\C^p$ with the base $\Omega$. Consider the class of
functions $f(z)$ on $T_\Omega$ that are approximated by sums
(\ref{sumexp}) with respect to the topology $\tau(T_\Omega)$ of
uniform convergence
on any domain $T_{D}$ with $D\Subset \Omega$. This class coincides with
the class $AP(T_\Omega)$ of holomorphic functions $f(z)$ in $T_\Omega$
such that the family $\{f(z+t)\}_{t\in\R^p}$ is a relatively compact
set in the topology $\tau(T_\Omega)$
\footnote{ It is easy to see that all sums (\ref{sumexp}) and
their limits in the topology $\tau(T_\Omega)$ belong to $AP(T_\Omega)$;
 for the converse assertion, see the Bochner--Fejer Theorem below.}.
 Such functions are called {\it holomorphic almost periodic on tube domains}.
We shall say that the {\it amoeba} $A_f$ of $\apf$ is the
closure of the projection
of the zero set of $f$ to $\Omega$. Note that the connected
components of the complement to the amoeba $A_f$ are the bases
$D_k\subset \Omega$ of the maximal tube domains
$T_{D_k}\subset T_\Omega$ without zeros of $f$;
every domain of this type is convex (see the book \cite{H1}, p.65).
Now the function $f(x+iy)$ is not periodic in each  variable
$x_1,\dots,x_p$, therefore we have to change the averages over
$[0,\,2\pi]$ of the increments in (\ref{ord1}) to their averages over
$[0,\,\infty]$.
In other words, define the order $\nu^k=(\nu_1^k,\dots,\nu_p^k)$
of any connected component $D_k$ of $\Omega\setminus A_f$
for $\apf$ by the formula
$$
\nu_j^k=\lim_{T\to \infty}{1\over 2T}\D_{|x_j|\le T}\Arg f(x+iy),
\quad y\in D_k, \quad \j1 p.
$$
(for simplicity we change $-D_k$ to $D_k$).

Holomorphic almost periodic functions in tube domains and their zero
sets were a subject of intensive study (see  \cite{FRR},
\cite{R0} -- \cite{R5}).
In particular, the notion of the Jessen function of an analytic almost
periodic function in a strip was extended to holomorphic almost
periodic functions in tube domains.
Namely, it was proved in \cite{R0} that the limit
\be
J_f(y)=\lim_{s\to\infty}\left({1\over 2s}\right)^p
\int_{|x_j|<s,\ \j1 p}\log|f(x+iy)|dx_1\dots dx_p
\label{Jes}
\en
exists for all $y\in D$ and $J_f(y)$ is a convex function; a relation
between $J_f(y)$ and the zero sets of $f$ was described, too. In the
paper \cite{R3} it was proved that the function $J_f(y)$ is linear on any
domain $D\Subset\Omega$ iff $f(z)\neq 0$ on $T_D$.
Moreover, in this domain the function $f$ has the form
\be
f(z)=\exp\{g(z)+i\<z,\,c(D)\>\}, \quad g\in AP(T_D), \quad
c(D)\in\R^p
\label{repr}
\en
(for $p=1$ see  \cite{JT}, p.188, for $p>1$ see
\cite{R3} and the Theorem below). Note that the
function $g(z)$ is uniformly bounded on $T(D')$ for every
$D'\Subset D$, therefore we have
\be
\nu^k=c(D_k)=-\grad J_f(y),\quad y\in D_k,
\label{grad}
\en
for any connected component $D_k$ of the set $\Omega\setminus A_f$.
The vector $c(D_k)$ is called the {\it mean motion of} $f$
in the domain $D_k$.

Now recall that the {\it spectrum} of $\apf$ is  the set
$$
\sp f=\{\lambda\in\R^p:\,
\lim_{s\to\infty}(2s)^{-p}
\int_{|x_j|<s,\j1 p}e^{-i\<z,\l\>}f(x+iy)\,dx\ne 0\},
$$
the spectrum is at most countable. Also, denote by $G(f)$
the minimal additive subgroup of $R^p$ containing $\sp f$, i.e.,
the set of all linear combinations of the elements of $\sp f$
with  integer coefficients.

The structure of amoebas for holomorphic almost periodic functions
is very complicated. For example, a sequence of different connected
components of $\OA$ can be condensed to an inner point of $\Omega$.
The problem of describing the mean motions for
almost periodic functions is very complicated, too. In the case $p=1$,
a complete  solution was given in the paper \cite{JT}, p. 229-259:

\medskip

{\bf Theorem JT.} {\sl For a convex function
$J(y)$ on $(a,\,b)$ to be the Jessen function of some $f\in AP(T_{(a,\,b)})$
with the spectrum in an additive countable subgroup $G\subset\R$,
it is sufficient and necessary that the following conditions be fulfilled:

i) for $y$ belonging to any interval of linearity of $J$,
 $J'(y)\in G$,

ii) for any $(\alpha,\,\beta)\Subset (a,\,b)$ there exist $k\in\N,\ K<\infty$,
and  numbers $\l^1,\dots,\l^k\in G$, linearly independent over $\Z$,
such that for arbitrary intervals of linearity $I_0,\,I_1
\subset(\alpha,\,\beta)$,
and points $y^0\in I_0,\ y^1\in I_1$, one has
\be
J'(y^1)-J'(y^0)=\sum_{j=1}^k r_j\l^j,\quad r=(r_1,\dots,r_k)\in\Q^k
\label{dif}
\en
and
$$
\|r\|\le K|J'(y^1)-J'(y^0)|,
\label{est}
$$

iii) if the group $G$ has a basis,
i.e., a generating set of elements linearly independent over $\Z$,
then the vector $r$ in (\ref{dif}) belongs to $\Z^k$ and
the number of distinct intervals of linearity intersecting
with $(\alpha,\,\beta)$ is finite.}

\medskip

For $p>1$ we know only two results of this kind. Namely,
in \cite{R3} the piecewise convex functions that are the Jessen functions of
holomorphic almost periodic functions were described; further, in \cite{R5}
it was proved that any holomorphic {\it periodic}
function $F$ of $q$ variables has a locally finite number of connected
components of the complement to its amoeba; this property remains true for
any almost periodic function which is the restriction of $F$
to a complex $p$-dimensional ($p<q$) hyperplane.

The aim of our paper is to prove the following result.

\medskip

{\bf Theorem.} {\sl For any $f\in AP(T_\Omega)$ the following
assertions are true:

i) for every connected component $D_0$ of $\OA$
the representation (\ref{repr}) is valid with $c(D_0)\in G(f)$,

ii) for every convex domain $D\Subset\Omega$ there exist $k\in\N$ and
$\l^1,\dots,\l^k\in G(f)$, linearly independent over $\Z$, such that
for any two connected components $D_0,\ D_1$ of $\DA$,
\be
c(D_1)-c(D_0)=\sum_{j=1}^k r_j\l^j,\quad r(D_1,\,D_0)=
(r_1,\dots,r_k)\in\Q^k
\label{dif1}
\en
and
\be
\|r(D_1,\,D_0)\|\le K\|c(D_1)-c(D_0)\|
\label{est1}
\en
with a constant $K<\infty$ depending only on $f$, $D$ and the
Lebesgue measure of the projection of the set $D_1-D_0$
to the unit sphere, i.e., the value
$\mes_{p-1}\{(y^1-y^0)/\|y^1-y^0\|:\,y^1\in D_1,\ y^0\in D_0\}$.
For a fixed component $D_0$, one can take
the same constant $K$ in (\ref{est1}) for all components $D_1$ of
$\DA$,

iii) moreover, if $G(f)$ has a basis, then the vector $r(D_1,\,D_0)$
in (\ref{dif1}) lies in $\Z^k$ and  the number of
connected components of $\DA$ is finite.}

\bigskip

For the reader's convenience, before prooving  the Theorem
we  present some known results which will be used in the proof.

{\bf 1.} (The Bochner--Fejer Approximation Theorem).
{\sl For every $\apf$ there exists  a sequence of exponential
polynomials (\ref{sumexp}) with $\l\in\sp f$ converging to $f$
in the topology $\tau(T_\Omega)$
(for $p=1$,  see \cite{C}, p. 79, or  \cite{JT}, p. 149;
in the case $p>1$ the proof is similar, see also \cite{R5}).}

{\bf 2.} (The Kronecker Theorem).
{\sl If the coordinates of  a vector $\mu\in\R^p$ are
linearly independent over $\Z$, then for any $a\in\R^p$ and
$\e>0$ there exist $t\in\R$ and $m\in\Z^p$ such that
\be
\|\mu t-a-2\pi m\|<\e.
\label{Kr}
\en
(see  \cite{B} or  \cite{C}, p. 147)}.

{\bf 3.} {\sl Any finitely generated subgroup of $\R^p$ has a finite basis
(see  \cite{H2}, p. 47)}.

\bigskip

{\bf Proof of the Theorem}.  Let $D_0$ be a connected component of
$\OA$ and $D'_0\Subset D_0$. Note that
\be
|f(z)|\ge\beta>0,\quad \forall z\in T_{D'_0}.
\label{nz3}
\en
Really,  otherwise
there exists a
sequence $z^n=x^n+iy^n\in T_{D'_0}$ such that
$f(z^n)\to 0$ as  $n\to\infty$.
Using the Bochner--Fejer Theorem and passing to a subsequence,
we get that the functions $f(z+x^n)$ converge uniformly
in $T(D_0)$ to a function
$g(z)$ and $y^n$ converge to a point $y'\in D_0$.
It is easy to see that $g(iy')=0$ and $g(z)\not\equiv 0$.
Consequently $f(z+x^n)$ has zeros in the set $T_{D_0}$
for $n$ large enough,  which is impossible.

It follows from (\ref{nz3}) and the Bochner-Fejer Theorem that
there exists an exponential polynomial $Q$ as in (\ref{sumexp}) such that
$\sp Q\subset\sp f$ and
\be
|f(z)-Q(z)|<|f(z)|/2,\quad z\in T_{D'_0}.
\label{apr1}
\en
Since $\sp Q$ is finite, we can choose a basis $\l^1,\dots,\l^q$
of the group $G(Q)$. Therefore we have
\be
Q(z)=\sum a_n \exp\{i\sum_{j=1}^q b_{n,j}\<z,\l^j\>\},
\quad a_n\in\C,\quad b_{n,j}\in\Z.
\label{pol0}
\en
Let $\oQ(\xi,y),\ \xi=(\xi_1,\dots,\xi_q)\in\R^q,\ y\in\R^p$,  be the function
\be
\oQ(\xi,y)=\sum a_n \exp\{-\sum_{j=1}^q b_{n,j}\<y,\l^j\>\}
 \exp\{i\sum_{j=1}^q b_{n,j}\xi_j\}.
\label{pol2}
\en
Denote by $\L$ the $q\times p$-matrix with the rows $\l^1,\dots,\l^q$.
If $x\in\R^p$ is not orthogonal to any linear combination
 of the vectors $\l^1,\dots,\l^q$ over $\Z$, then the components
of the vector $\L x\in\R^p$ are linearly independent over $\Z$ and,
by the Kronecker Theorem, the set
$$
\{\L x+2\pi m:\,x\in\R^p,\,m\in\Z^q\}
$$
is dense in $\R^q$. Further,
the function $\oQ(\xi,y)$ has the period $2\pi$
in each of the variables $\xi_1,\dots,\xi_q$. Using (\ref{apr1}),
(\ref{nz3}), and the identity $Q(z)=\oQ(\L x,y)$, we obtain
\be
|\oQ(\xi,y)|\ge\beta/2>0
\quad\forall \xi\in\R^q,\  y\in D'_0.
\label{nz4}
\en
Denote by $\D$ the vector in $\R^q$ whose components
$\D_j,\ \j1 q$, equal the increments of $\Arg\oQ(\xi,y)$ as $\xi_j$ runs
over the segment $[0,\,2\pi]$, $\xi_l$
fixed for $l\neq j$, $y\in D'_0$.
The function $\oQ(\xi,y)$ is periodic, hence $\D_j\in 2\pi\Z$
for all $j$. Taking into account the continuity of $\oQ$,
we see that these increments depend  on neither $\xi_l$ nor $y\in D'_0$.
In view of (\ref{nz4}), we can define the continuous function
$$
\overline{h}(\xi,y)=\log\oQ(\xi,y)
-{i\over 2\pi}\sum_{j=1}^q \D_j\xi_j
+{1\over 2\pi}\sum_{j=1}^q\D_j\<y,\l^j\>.
$$
This function has the period $2\pi$ in each of  the variables $\xi_j$,
hence it is bounded on $\R^q\times D'_0$.
It is obvious that for $h(z)=\overline{h}(\L x,y)$ we have
\be
\exp h(z)=Q(z)\exp\{-i\<z,\,\L'{\D\over 2\pi}\>\},
\label{exph}
\en
where $\L'$ is the transpose of $\L$. We see that the function
(\ref{exph}) belongs to $AP(T_D)$ and the function $h(z)$ is
bounded on $T_{D'_0}$. Now, using the definition of holomorphic almost
 periodic function, we obtain $h(z)\in AP(T_{D'_0})$.
Further, using (\ref{apr1}), we get that
the function $\hat h=\log[f(z)/Q(z)]$ is well-defined,
almost periodic and bounded in $T_{D'_0}$. Therefore, we have
\be
f(z)=\exp\{\hat h(z)+ h(z)+i\<z,\,\L'{\D\over 2\pi}\>\}.
\label{exph1}
\en
Using (\ref{Jes}), (\ref{grad}) for the domain $D_0$, and the boundedness
of the functions $\hat h(z)$ and $h(z)$ in $T_{D'_0}$, we
prove assertion i) with $c(D_0)=\L'\D/2\pi$.

To prove ii), first suppose that the convex domain
$D'\Subset\Omega$ has  the property
\be
\forall \yD\quad |f(x^0+iy)|\ge 5\gamma>0,
\label{nz1}
\en
with some $x^0\in\R^p$; we may assume that $x^0=0$.
Using the Bochner-Fejer Theorem, we can take  an exponential polynomial
$P$ as in (\ref{sumexp}) such that $\sp P\subset\sp f$ and
\be
|f(z)-P(z)|<\gamma,\quad\forall z\in T_{D'}.
\label{apr}
\en
Since $\sp P$ is finite, we can choose a basis $\l^1,\dots,\l^k$
of the group $G(P)$. Therefore we have
$$
P(z)=\sum a'_n \exp\{i\sum_{j=1}^k b'_{n,j}\<z,\l^j\>\},
\quad a'_n\in\C,\quad b'_{n,j}\in\Z.
$$
Hence, for any $z\in T_{D'},\ x\in\R^p$ we clearly have
\be
|P(z+x)-P(z)|\le\sum |a'_n| \sup_{\yD}\exp\{-\sum_{j=1}^k b'_{n,j}\<y,\l^j\>\}
|1-\exp\{i\sum_{j=1}^k b'_{n,j}\<x,\l^j\>\}|.
\label{ine1}
\en

Let $\tL$ be the $k\times p$-matrix with the rows $\l^1,\dots,\l^k$,
let $\e>0$, and
\be
F=\{x\in\R^p:\,\exists \tm\in\Z^k\ \|\tL x-2\pi \tm\|<3\e\}.
\label{F}
\en
It follows from (\ref{ine1}) that for $\e$ small enough and $x\in F$,
\be
|P(z+x)-P(z)|<\gamma.
\label{ine2}
\en
Combining (\ref{nz1}) with $x^0=0$, (\ref{apr}), and (\ref{ine2}),
we obtain
\be
|f(x+iy)|\ge 2\gamma>0 \quad \forall x\in F,\quad \yD.
\label{nz2}
\en

Let $D_0,\ D_1$ be any distinct connected components of the set
$D'\setminus A_f$,  and let $y^0\in D_0,\ y^1\in D_1,\
\nu\in\R^p,\,\|\nu\|=1$ be such that for some $s>0$
\be
y^1=y^0+\nu s
\label{points}
\en
and, in addition,
\be
\<\l,\nu\>\neq 0 \quad \forall \l\in G(f).
\label{nort}
\en
As before, we have
\be
|f(x+iy^l)|\ge\beta'>0\quad\forall x\in\R^p,\quad  l=1,2.
\label{nz5}
\en

It follows from (\ref{nz3}) and the Bochner-Fejer Theorem that
there exists an exponential polynomial $Q$ such that
$\sp Q\subset\sp f$,
\be
|f(z)-Q(z)|<\gamma,\quad z\in T_D',
\label{apr2}
\en
and the inequality (\ref{apr1}) is fulfilled for all
$x\in\R^p,\ y=y^l,\ l=1,2.$

Note that the polynomial $Q$ depends on the domains $D_0,\ D_1$,
while the polynomial $P$  and vectors $\l^1,\dots,\l^k$ depend on the
domain $D'$ and don't depend on $D_0$ and $D_1$.

Let the dimension of the linear span  of $\sp P\cup\sp Q$ over $\Q$
be $q\ge k$ and let $\l^{k+1},\dots,\l^q$ be vectors that together with
the vectors $\l^1,\dots,\l^k$ form a basis of this span. Choose
$N\in\N$ such that any element of $\sp Q$ is a linear combination
 of $\l^1/N,\dots,\l^k/N,\,\l^{k+1},\dots,\l^q$ over $\Z$. If the group
$G(f)$ has a basis, then we can take at first $\l^1,\dots,\l^k$, next
$\l^{k+1},\dots,\l^q$ from this basis and take always $N=1$.

Let $\L$ be the $q\times p$-matrix with the rows $\l^1,\dots,\l^q$.
It is obvious that $Q(z)$ can be written in the form (\ref{pol0}) with
$b_{n,j}/N$ instead of $b_{n,j}$. Let $\oQ(\xi,y)$ be defined
by formula (\ref{pol2}) with the same changes. As  above, we have
$$
Q(z)=\oQ(\L x,y),
$$
but now the function $\oQ(\xi,y)$ has the period $2\pi N$
in each of the variables $\xi_1,\dots,\xi_q$.
Using the Kronecker Theorem with $\mu=\L x/N,\ a=\xi/N$,
we find that the set
$$
\{\L x+2\pi Nm:\,x\in\R^p,\,m\in\Z^q\}
$$
is dense in $\R^q$. As before, it follows from (\ref{nz5}) that
\be
|\oQ(\xi,y^l)|\ge\beta'/2>0\quad \forall\xi\in\R^q,\ l=1,2.
\label{nz6}
\en

Denote by $\Pr[\xi]$ the vector formed  by the first $k$ coordinates of
$\xi\in\R^q$. Let
$$
E=\{\xi\in\R^q:\,\exists m\in\Z^q \quad\|\Pr[\xi-2\pi m]\|<3\e\},
$$
where $\e$ is the same as in (\ref{F}). From the Kronecker Theorem
it follows that the set $\{\L x+2\pi Nm:\,x\in\R^p,\,m\in\Z^q\}$
is dense in $E$. Since $F=\{x\in\R^p:\,\L x\in E\}$,
we see that the set $\{\L x+2\pi Nm:\,x\in F,\,m\in\Z^q\}$
is dense in $E$, too.
Therefore, taking into account (\ref{nz2})
 and (\ref{apr2}), we obtain
\be
|\oQ(\xi,y)|\ge\gamma>0\quad\forall\xi\in E,\quad \yD.
\label{nz0}
\en

Denote by $\D^\rho$ the vectors in $\R^q$ whose components
$\D_j^\rho,\ \j1 q$, equal the increments of the arguments for the functions
$\oQ(\xi,y^\rho)$ while $\xi_j$ run the segment $[0,\,2\pi N]$, $\xi_l$
fixed for $l\neq j$,  $\rho=0,1$. Arguing as above, we define
the bounded periodic functions
$$
\overline{h}_\rho(\xi)=\log\oQ(\xi,y^\rho)
-{i\over 2\pi}\sum_{j=1}^q \D^\rho_j{\xi_j\over N}
+{1\over 2\pi N}\sum_{j=1}^q\D^\rho_j\<y^\rho,\l^j\>,\ \rho=1,2.
$$
Using (\ref{apr1}) with $z=x+iy^\rho$, we obtain that
equality (\ref{exph1}) is fulfilled with
$h(z)=\overline{h}_\rho(\L x),\ \hat h(z)=\log[f(z)/Q(z)]$,
and $z=x+iy^\rho,\ \rho=0,1$. Thus we have
\be
c(D_\rho)={1\over 2\pi N}\L'\D^\rho,\quad \rho=0,1.
\label{mm}
\en

Further, by $\d(\xi^0,\xi^1)$ denote the increment of  the
argument of the function $\oQ(\xi,y)$ over the boundary of the rectangle
$\Pi(\xi^0,\xi^1)$ with vertices at
$\xi^0+iy^1,\ \xi^1+iy^1,\ \xi^1+iy^0,\ \xi^0+iy^0$,
which run in the indicated order. It follows from (\ref{nz0}) and
(\ref{nz6}) that $\d(\xi^1,\xi^2)$ is defined for arbitrary
points $\xi^0,\ \xi^1\in E$. We obviously have
\be
\d(\xi^0,\xi^2)=\d(\xi^0,\xi^1)+\d(\xi^1,\xi^2)
\quad\forall\xi^0,\xi^1,\xi^2\in E.
\label{tr}
\en

Moreover, it follows from the periodicity of $\oQ(\xi,y^\rho)$
that in order to compute $\d(0,\xi)$ with $\xi\in 2\pi N\Z^q$
we only need to take into account the increments in the planes
$y=y^0$ and $y=y^1$. Hence for every $m=(m_1,\dots,m_q)$ we obtain
\be
\d(0,2\pi Nm)=\sum_{j=1}^q m_j\D_j^1-\sum_{j=1}^q m_j\D_j^0=
\<m,\D^1-\D^0\>.
\label{d1}
\en
Further, if $\xi'\in\R^q$ and $m\in\Z^q$ satisfy
$\|\Pr[\xi'-2\pi m]\|<3\e$, then $\xi'\in E$ and the rectangle
$\Pi(\xi',2\pi m)$ can be contracted inside the set $E\times[y^0,y^1]$ to
the segment $2\pi m\times[y^0,y^1]$. Therefore, we have
\be
\d(\xi',2\pi m)=0.
\label{d2}
\en
In particular, if $\xi',\xi''\in E$ and $\Pr[\xi'-\xi'']=0$,
then $\d(\xi',\xi'')=0$. Hence $\D_j^1=\D_j^0$ for $k<j\le q$
and
\be
\D^1-\D^0=(2\pi d,\ 0), \quad d\in\Z^k,\quad 0\in\Z^{q-k}.
\label{d}
\en

Moreover, it follows from (\ref{pol2})  that the function
$Q(\xi+\L\nu u,y^0+\nu v)$ is analytic  with respect to the variable
$w=u+iv$. From (\ref{points}) and (\ref{apr1}) it follows that this
function has no zeros
for $v=0$ and $v=s$. By the Argument Principle,
for each $\xi\in E$ and $u>0$ such that $\xi+\L\nu u\in E$, we have
\be
\d(\xi,\xi+\L\nu u)\ge 0.
\label{d3}
\en

Note that the components $\mu_j$ of the vector
$\tL\nu$ are linearly independent over $\Z$, hence
inequality (\ref{Kr}) with $\mu=(\mu_1,\dots,\mu_k)$ has a solution
for every $a=(a_1,\dots,a_k)\in\R^k$ and $\e>0$.
This means that the almost periodic function
$Y(t)=e^{i(\mu_1t-a_1)}+\dots+e^{i(\mu_kt-a_k)}$
takes values arbitrarily close to $k$. It follows from Bohr's definition
of almost periodic functions (see the book \cite{C}, p. 14)
 that for every $\eta>0$ the inequality $|Y(t)-k|<\eta$ has a solution
on every interval of fixed length $M$. Hence  inequality (\ref{Kr})
has a positive solution and each point $a\in\R^k$ belongs
to the union of the open sets
$$
H_l=\{a\in\R^k:\,\exists t'\in(0,\,l),\
\exists\tm\in\Z^k,\  \|\tL\nu t'-a-2\pi \tm\|<\e\}.
$$
Therefore  $H_L\supset [0,\,2\pi]^k$ for some $L<\infty$ and
inequality (\ref{Kr}) with $\mu=\tL\nu$ has a solution on the segment
$[0,\,L]$ for all $a\in\R^k$.

Take $t\in(0,\,L)$ and $\tm^0\in\Z^k$ such that
\be
\|\tL\nu t-2\e {d\over \|d\|}-2\pi \tm^0\|<\e.
\label{ine3}
\en
 Observe that this inequality implies
$\|\tL\nu t-2\pi \tm^0\|<3\e$. It follows from (\ref{d2}) that
for any $m\in\Z^q$ and $m^0=(\tm^0,\,0)\in\Z^q$,
\be
\d(2\pi m+\L\nu t,2\pi m+2\pi m^0)=0.
\label{d4}
\en
By (\ref{d3}), we have
\be
\d(2\pi m, 2\pi m+\L\nu t)\ge 0.
\label{d5}
\en
Combining (\ref{tr}), (\ref{d4}), and (\ref{d5}), we get
\be
\d(0,2\pi Nm^0)=\sum_{n=0}^{N-1}[\d(2\pi nm^0, 2\pi nm^0+\L\nu t)+
\d(2\pi nm^0+\L\nu t,2\pi (n+1)m^0)]\ge 0.
\label{d6}
\en
Using (\ref{d6}), (\ref{d}), and (\ref{d1}), we have
\be
2\pi\<\tm^0,d\>=\<m^0,\D^1-\D^0\>=\d(0,2\pi Nm^0)\ge 0.
\label{d7}
\en
 Then, by (\ref{ine3}) and (\ref{d7}),
\be
\<\tL\nu t,d\>=\<\tL\nu t-{2\e d\over \|d\|}-2\pi \tm^0,d\>+
2\pi \<\tm^0,d\>+2\e\|d\|\ge \e\|d\|.
\label{d8}
\en
Taking into account (\ref{mm}) and (\ref{d}), we obtain
\be
c(D_1)-c(D_0)={1\over 2\pi N}\L'(\D^1-\D^0)=\tL'd/N.
\label{mm2}
\en
Hence we have proved (\ref{dif1}) with
\be
r(D_1,D_0)=d/N.
\label{mm3}
\en
Combining (\ref{d8}), (\ref{mm2}), (\ref{mm3}) with the equality
$\|\nu\|=1$, we get
$$
\|r(D_1,D_0)\|=\|d\|/N\le {t\over \e N}|\<\tL\nu,d\>|\le
{L\over\e N}|\<\nu,\tL'd\>|\le (L/\e)\|c(D_1)-c(D_0)\|.
\label{est3}
$$

Thus we  have proved assertion ii) for domains $D\Subset \Omega$
 satisfying (\ref{nz1}) and connected components
$D_0,\ D_1$ of $\DA$. Note that $k$ and $\l^1,\dots,\l^k$ in (\ref{dif1})
depend on $D$ only, and $K$ in (\ref{est1}) depends on $D$ and $\nu$.

Now let $\l^1, \l^2, \dots$, be a basis of the linear span of $\sp f$
over $\Q$. It is clear that all  $\l^j$ can be taken  from $G(f)$.
It follows from assertion i) that for each connected component $\tilde D$
of $\DA$ there exist $k=k(\tilde D)\in\N$ such that
\be
c(\tilde D)=\sum_{j=1}^kr_j\l^j,\quad r_j\in\Q,\quad\forall j,
\quad r_k\neq 0.
\label{repr1}
\en
Let us show that $k(\tilde D)$ are uniformly bounded for all $\tilde D$.
Assume the contrary. Then there exists a sequence of connected
components $D_n$ of $\DA$ such that $k(D_n)\to\infty$ as $n\to\infty$.
Take $y^n\in D_n,\ n=1,2,\dots$.
We may assume that
$y^n\to y'\in\overline{D}$. Since $f(z)\not\equiv 0$ on $T_D$,
we obtain $f(x+iy')\not\equiv 0$ for $x\in\R^p$. Hence
we can take $x^0,\ \gamma>0$ and a convex neighborhood $D'$ of $y'$
such that (\ref{nz1})  is fulfilled. Then there exists $n_0$ such that
$D_n\cap D'\neq\emptyset$ for all $n\ge n_0$. Using (\ref{dif1}) for
the domains $D_n\cap D',\ D_{n_0}\cap D'$, we get
$$
c(D_n)-c(D_{n_0})=\sum_{j=1}^{k'}r_j\l'^j,\quad
(r_1,\dots,r_{k'})=r(D_n,D_{n_0})\in\Q^{k'}
$$
where $\l'^j,\,\j1 k'$ belong to the linear span of $\sp f$
and the number $k'$ is the same for all $n$.
This contradicts the unboundedness of $k$ in (\ref{repr1}).

Now it follows from (\ref{repr1}) that (\ref{dif1}) is valid for all
components $D_1,\ D_0$ of $\DA$ with $k$ depending on $D$ only.

Furthermore, assume that (\ref{est1}) is false. This means that
there exist two sequences $\{D_n\},\ \{D'_n\}$  of connected components
of $\DA$ such that
\be
\| c(D_n)-c(D'_n)\|/\|r(D_n,D'_n)\|\to 0
\label{to}
\en
as $n\to\infty$ and
$$
{\rm mes}_{p-1}\kappa(D_n-D'_n)\ge\beta>0,
$$
where $\kappa(y)=y/\|y\|.$
Since (\ref{nort}) is false only for a countable number
of hyperplanes in $\R^p$, we see that there exists
$\nu\in\bigcap_{l=1}^\infty\bigcup_{n=l}^\infty
\kappa(D_n-D'_n)$ satisfying (\ref{nort}). Take a subsequence
$n_l$ such that $\nu\in\kappa(D_{n_l}-D'_{n_l}), \ l=1,2,\dots$.
Omit the first indices. Without loss of generality
it can be assumed that there exist $y^l\in D_l,\
y'^l\in D'_l$, and $s_l>0$ such that
$y^l\to \tilde y\in\overline{D},\ y'^l\to y'\in\overline{D}$
and $y^l= y'^l+\nu s_l$. Clearly, we have
$\tilde y= y'+\nu s, \ s\ge 0$.

Further, since $f(x+iy')\not\equiv 0$ for $x\in\R^p$, we can take
$x^0\in\R^p$ such that $f(x^0+iy')\neq 0$. Hence the
function $g(w)= f(x^0+\nu w+iy')$, analytic
in a neighborhood of set $\{w=u+iv:\,0\le v\le s,\ u\in\R\}$,
 is not identically zero.
Therefore $g(u_0+iv)\neq 0$ for some $u_0\in\R$ and all
$v\in[0,\,s]$. Consequently, the function
$f(x_0+\nu u_0+iy)$ does not vanish on the segment
$[y',\,\tilde y]\in\overline{D}$ (if $s=0$, then this segment is
a point and we can take $u_0=0$).

Take a convex domain $D'$, with $[y',\,\tilde y]\subset D'\Subset\Omega$
such that (\ref{nz1}) is true for this domain with $x^0+\nu u_0$
instead of $x^0$. The domain $D'$ intersects the domains
$D_l,\ D'_l$ for $l$ large enough. Now we can apply
(\ref{dif1}) and (\ref{est1}) to the domains $D_l\cap D',\ D'_l\cap D'$
and  the vector $\nu$. We have
\be
c(D_l)-c(D'_l)=\sum_{j=1}^{k'}r'_j\l'^j,\quad
(r'_1,\dots,r'_{k'})=r'(D_l,D'_l)\in\Q^{k'}
\label{dif3}
\en
and
\be
\|r'(D_l,D'_l)\|\le K'\|c(D_l)-c(D'_l)\|,
\label{est2}
\en
where $K'$ is the same for all $n$.
Pass in (\ref{dif3}) to the previous basis $\l^1,\l^2,\dots$.
Taking into account the obvious inequality
$$
\|r(D_l,D'_l)\|\le C(\l^1,\l^2,\dots,\l'^1,\l'^2,\dots) \|r'(D_l,D'_l)\|,
$$
we obtain that (\ref{est2}) contradicts (\ref{to}).
This  proves the first part of ii).

Finally note that for every component $D_0$ of $\DA$ we have
$$
\mes_{p-1}\kappa(y-D_0)\ge\eta(D_0)>0, \quad\forall y\in D.
$$
Hence for  a fixed $D_0$, the constant $K$ in (\ref{est1}) does not
depend on $D_1$. This concludes the proof of assertion ii).

As  was mentioned above, for the group $G(f)$ with a basis
we can take $\l^1,\dots,\l^k$ in (\ref{repr1}) from this basis and, moreover,
$N=1$ in (\ref{mm}). It follows from (\ref{mm3}) that the components
of  the vectors $r(D_0,\,D_1)$ are integers for all domains $D_0, \,D_1$.
Taking into account (\ref{grad}) and the convexity of $J_f(y)$, we obtain
that  the mean motions $c(D_k)$ are uniformly bounded for all components
$D_k$ of $\DA$. Using inequality (\ref{est1})  with a fixed $D_0$, we see
that $c(D_k)$ takes only a finite number of values.
Further, it follows from the convexity of $J_f$
that $c(D_k)\neq c(D_j)$ for any distinct components.
Hence there exists only a finite number of components $D_k$.
The theorem is proved.

\medskip

{\bf Remark.} Actually we have proved the following result:

{\sl For a convex function $J(y)$ on a convex domain $\Omega$ to be the
Jessen function of some $f\in AP(T_{\Omega})$
with the spectrum in an additive countable subgroup $G\subset\R^m$,
it is necessary that the following conditions be fulfilled:

i) if $D_0$ is the domain of linearity of $J$ and
$c(D_0)=-\grad J(y),\quad y\in D_0$, then $c(D_0)\in G$,

ii) if the group $G$ has a basis,
then the number of distinct domains of linearity intersecting
with any fix domain $D\Subset\Omega$ is finite;

iii) in general case, for any convex domain $D\Subset\Omega$ there exist
$k\in\N$ and vectors $\l^1,\dots,\l^k\in G$
such that for any two distinct domains $D_0,\ D_1\subset D$ of linearity
of $J$
\be
c(D_1)-c(D_0)=\sum_{j=1}^k r_j\l^j,\quad r=r(D_1,D_0)=(r_1,\dots,r_k)\in\Q^k;
\en
moreover, for each $\nu\in\{y\in\R^m:\,\|y\|=1\}\setminus
\cup_{\l\in G}\{y:\,\<y,\l\>=0\}$ there exists a constant $K(\nu)<\infty$
with the property}
$$
\nu\in\{(y^1-y^0)/\|y^1-y^0\|:\,y^1\in D_1,\ y^0\in D_0\}\Rightarrow
\|r(D_1,D_0)\|\le K(\nu)\|c(D_1)-c(D_0)\|.
$$

I don't know if the word "necessary" can be changed to the words
"necessary and sufficient" here.

\medskip

{\bf Acknowledgements.} The author is very grateful to
Professor A.Yu.\,Rashkovskii for useful discussions.

\bigskip
{\it Kharkiv National University

Department of Mechanics and Mathematics

md. Svobody 4, Kharkiv 61077, Ukraine

\medskip
e-mail: favorov@ilt.kharkov.ua}

\vskip 1.5cm


\bigskip

\begin{thebibliography}{}

\bibitem{B} {\bf H.\,Bohr.} Again the Kronecker Theorem.
Journ. of the London Math. Soc. {\bf 9} (1934), No.33-36, p.5.

\bibitem{C} {\bf C.\,Corduneanu.} Almost periodic functions.
Interscience Publishers, New-York -- London -- Sydney -- Toronto,
a division of John Wiley, 1961.

\bibitem{FRR}
{\bf S.Yu.\,Favorov, A.Yu.\,Rashkovskii and L.I.\,Ronkin.} Almost
periodic currents and holomorphic chains. C. R. Acad.
Sci. Paris {\bf 327}, Serie I (1998), 302-307.


\bibitem{P}
{\bf M.\,Forsberg, M.\,Passare, A.\,Tsikh}. Laurent determinants
and arrangement of hyperplane amoebas. Adv. Math. {\bf 151} (2000),
no. 1, 45-70.

\bibitem{M}
{\bf G.\,Mikhalkin}.
Real algebraic curves, the moment map and amoebas.
[J] Ann. Math. (2) {\bf 151} (2000), no.1, 309-326.


\bibitem{G} {\bf I.\,Gelfand, M.\,Kapranov, and A.\,Zelevinsky}.
Discriminats, resultants and multidimensional determinants.
Boston, Birkh\"auser, 1994.

\bibitem{H2}
{\bf M.\,Hall}.  The theory of groups. New York, The Macmillan
Company, 1959.

\bibitem{H1}
{\bf L.\,H\"ormander.}  An Introduction to Complex Analysis in Several
Variables. D.~van Nostrand company,  New Jersey, Inc. Princeton, 1966.


\bibitem {JT}
{\bf B.\,Jessen and H.\,Tornehave}. Mean motion and zeros of
almost-periodic functions. Acta Math. {\bf 77} (1945), 137-279.


\bibitem{R0}
{\bf L.I.\,Ronkin.} Jessen's theorem for holomorphic almost periodic
functions in tube domains. Sibirsk. Mat. Zh. {\bf 28} (1987), 199-204
(Russian).

\bibitem{R1}
{\bf L.I.\,Ronkin.} Jessen's theorem for holomorphic almost
periodic mappings.
Ukrainsk. Mat. Zh. {\bf 42} (1990), 1094-1107 (Russian).

\bibitem{R2}
{\bf L.I.\,Ronkin} Functions of Completely Regular Growth.
Mathematics and its Appl., Soviet Ser., v. 81.
Kluwer Acad. Publ., Dordrecht Boston, 1992.

\bibitem{R3}
{\bf L.I.\,Ronkin.} On a certain class of holomorphic almost periodic
functions.
Sibirskii Mat. Zh. {\bf 33} (1992), 135-141 (Russian).

\bibitem{R4}
{\bf L.I.\,Ronkin.} Almost periodic distributions and divisors
in tube domains. Zap. Nauchn. Sem. POMI {\bf 247} (1997), 210-236 (Russian).

\bibitem{R5}
{\bf L.I.\,Ronkin.} On zeros of almost periodic functions generated
by holomorphic functions in a multicircular domain. In: Complex
Analysis in Modern Mathematics, Moscow, Fazis (2001), 243-256.

\bibitem{Ru}
{\bf H.\,Rullga ard.} Stratification des espaces de polyn\^omes de
Laurent et la structure de leurs amibes. C. R. Acad. Sci. Paris
{\bf 331}, Serie I (2000), 355-358.

\end{thebibliography}
\end{document}